\documentclass[a4paper,11pt]{article}
\usepackage{amsfonts}
\usepackage{amssymb}
\usepackage{amsmath}
\usepackage{amsthm}
\usepackage{multirow}
\usepackage{diagbox}
\usepackage{amsopn}
\usepackage{gensymb} 
\usepackage{fancyhdr}
\usepackage{graphicx}
\usepackage{mathrsfs}
\usepackage{latexsym}
\usepackage{graphicx}
\usepackage{subfigure}
\usepackage{color}
 \usepackage{calc}
  \usepackage{tikz}
  
\usepackage{setspace}
\usepackage{textcomp} 

\usepackage[width=11cm]{geometry}   
\usepackage{varwidth}
\usepackage{algorithm}
\usepackage{multicol}
\usepackage{algpseudocode}
\numberwithin{equation}{section} \allowdisplaybreaks
\theoremstyle{plain}
\newtheorem{thm}{Theorem}[section]
 \usepackage{rotating}
\usepackage{appendix}

\theoremstyle{remark}
\newtheorem{remark}{Remark}
\numberwithin{equation}{section}
\newtheorem{Lemma}[thm]{Lemma}

\newtheorem*{Proof}{Proof}

\begin{document}
\title{\bf{Bayesian approach for inverse obstacle scattering with Poisson data}}

\author{\hspace{-2.5cm}{\small Xiao-Mei Yang$^{1}$ and Zhi-Liang Deng$^2$\thanks{
Corresponding author: dengzhl@uestc.edu.cn
Supported by NSFC No. 11601067, 11771068 and 11501087, the Fundamental Research Fund for the Central Universities No. ZYGX2018J085.
}} \\
\hspace{-1.5cm}{\scriptsize $1.$  School of Mathematics,
Southwest Jiaotong University,
Chengdu 610031, China}\\
\hspace{-2.5cm}{\scriptsize $2.$ School of Mathematical Sciences,  University of Electronic Science and Technology of China,
Chengdu 610054, China}
}
\date{}
\maketitle

\begin{abstract}
\noindent We consider an acoustic obstacle reconstruction problem with Poisson data. Due to the stochastic nature of the data, 
we tackle this problem in the framework of Bayesian inversion. 
The unknown obstacle is parameterized in its angular form. The prior for the parameterized unknown plays key role in the Bayes reconstruction algorithm. The most popular used prior is the Gaussian. Under the Gaussian prior assumption, we further suppose that the unknown satisfies the total variation prior. With the hybrid prior, the well-posedness of the posterior distribution is discussed.
The numerical examples verify the effectiveness of the proposed algorithm. 




\end{abstract}

\section{Introduction}
Inverse obstacle scattering problems occur in many real-world applications such as radar, sonar and non-destructive testing \cite{colton}.  One main object of this kind of problems is to find the shape and location of the unknown obstacle, which usually can not be observed directly. Due to the widely applications, these prolems have attracted many researchers' attention. Especially, a large number of papers put emphasis on designing some effect numerical algorithms  to recover the unknown information according to some related observed data \cite{bui, colton}.
The observed data are contaminated by the noise inevitably.  The Gaussian process has many theoretical and practical advantages since it can mimic the effect of many random processes. And therefore a lot of noise models can be characterized by the Gaussians.  
  However, in many applications, the response variable of interest is a count, that is, takes on nonnegative integer values. For count data, it is well modeled by a Poisson process. 
The Poisson process is usually used to describe the occurrences of events and the time points at which the events occur in a given time interval, such as the occurrence of natural disasters and the arrival times of customers at a service center.
In the obstacle scattering problem, the counts of photons at a remote closed surface are detected. As seen in \cite{hohage,werner, zhou}, the counts of photons have interacted with the unknown object and can be well modeled by a Poisson process. In view of this, we consider to use the Poisson data to recover the unknown shape. 

The stochastic nature of the data undoubtedly results in the uncertainty of the reconstruction. Therefore, besides the reconstruction of the unknown object itself,  its uncertainty requires the extra concerns. Probabilistic thinking provides a natural way to quantify the uncertainty. Recently, Bayesian inference method has become a popular tool for this purpose \cite{stuart1}. In Bayesian approach, the main idea is to translate the prior knowledge on the noise and on the unknowns to prior probability laws and then use the forward model to transmit the prior information to the posterior distribution from which we can deduce any
information about the unknowns. A rigorous Bayesian framework is developed for the inverse problems with the infinite dimensional unknowns \cite{stuart2}. From the existed works, we can see that the prior for the unknowns is crucial in Bayesian approach and even determines whether or not it works. A difficult subject is how to code the experience prior information in some distribution formula.  
 It is well-known that, subject to certain mild conditions,  the sum of a set of random variables, has a distribution that becomes increasingly Gaussian as the number of  terms in the sum increases. 
 This point allows us to use Gaussians to model the prior for many real physical parameters.
In this paper, we characterize the unknown by two kinds of prior knowledges: the Gaussian and the hybrid priors. Upon these priors, we discuss the well-posedness of the posterior distribution with the Poisson data. 
 
The remainder of this paper is organized in the following: In Section 2, we describe the exterior acoustic scattering problem. In Section 3, the Bayesian approach and the well-posedness are discussed. In the final, some numerical examples are given. 

%

\section{Acoustic obstacle scattering with Poisson data}
For a perfect cylindrical conductor with smooth cross section $D\subset\mathbb{R}^2$, the scattering of polarized, transverse magnetic time harmonic electromagnetic waves is described by the Helmholtz equation
\begin{align}
\triangle u+k^2u=0 \,\, \text{in} \,\, \mathbb{R}^2\backslash D,  \label{ac2.1}
\end{align}
where $k>0$ is the wave number, $u=u^{\rm s}+u^{\rm i}$  the total field, $u^{\rm s}$ the scattered field and $u^{\rm i}$ the incident wave. We let $u^{\rm i}=\exp(ikx\cdot d)$, where $d$ is the incident direction. The following boundary conditions are imposed
\begin{align}
&u\mid_{\partial D}=0, \label{ac2.2}\\
&\lim_{r\rightarrow \infty}\sqrt{r}(\frac{\partial u^{\rm s}}{\partial r}-iku^{\rm s})=0,\,\, r=|x|. \label{ac2.3}
\end{align}
The condition \eqref{ac2.2} corresponds to the sound soft boundary condtion and  \eqref{ac2.3} is the so-called Sommerfeld radiation condition.  
The Sommerfeld radiation condition characterizes the outgoing wave and enables us to have
the following asymptotic behaviour for the scattered wave $u^{\rm s}$:
\begin{align}\label{ac2.4}
u^{\rm s}(x)=\frac{\mathrm{e}^{i k|x|}}{|x|^{\frac{1}{2}}}\left\{u_{\infty}(\hat{x})+\mathcal{O}\left(\frac{1}{|x|}\right)\right\} \quad \text{as} \quad|x| \rightarrow \infty,
\end{align}
where $\hat{x}=\frac{x}{|x|}$ and $u_{\infty}(\hat{x})$ is the far field pattern. 
The direct scattering problem is to find the scattered field $u^{\rm s}$ for given incident field $u^{\rm i}$ and obstacle $D$. It is known  that there exists a unique solution $u\in H^1(\mathbb{R}^2\backslash D)$ if $\partial D$ is Lipschitz continuous \cite{colton}.

We focus on the reconstruction of the unknown domain $D$ by some observed data:  the photon counts of the scattered electromagnetic field far away from the obstacle. Since at large distances the photon density is approximately proportional to $|u_\infty|^2$ \cite{werner}, the obstacle reconstruction problem can be described by the operator equation
\begin{align}\label{ac2.5}
F(D)=\tau |u_\infty|^2,
\end{align}
where $\tau$ is a constant corresponding the noise level. 
We suppose that  the domain $D$ is starlike, i.e., the boundary $\partial D$ is parameterized by
\begin{align}
\label{ac2.6}
\partial D=q(t)[\cos(t), \sin(t)]^T, \quad t\in[0, 2\pi).
\end{align}
In this setting, the function $q$ becomes the unknown and therefore the equation  \eqref{ac2.5} is equivalent to 
\begin{align}\label{ac2.7}
G(q)=\tau |u_\infty|^2. 
\end{align}
The operator equation can be given by the single-double layer potential theories \cite{colton}. Let $\partial D$ be of class $C^2$. Define the single layer potential operator $S$
\begin{align}\label{ac2.8}
(S \varphi)(x)=2 \int_{\partial D} \Phi(x, \tilde{x}) \varphi(\tilde{x}) d s(\tilde{x}), \quad x \in \partial D,
\end{align}
and the double layer potential operator $K$
\begin{align}\label{ac2.9}
(K \varphi)(x)=2 \int_{\partial D} \frac{\partial \Phi(x, \tilde{x})}{\partial \nu(\tilde{x})} \varphi(\tilde{x}) d s(\tilde{x}), \quad x \in \partial D,
\end{align}
where $\Phi(x, \tilde{x})=\frac{i}{4} H_{0}^{1}(|x-\tilde{x}|)$ is the fundamental solution, $H_{0}^{1}(\cdot)$ is the Hankel function of the first kind of order zero and $\nu(\tilde{x})$ is normal direction. It can be seen in \cite{colton} that $S$ and $K$ are bounded from $C^{0, \alpha}(\partial D)$ into $C^{1, \alpha}(\partial D)$. With the single-double layer potentials, the scattered field can be written as
\begin{align}\label{ac2.10}
u^{\rm s}(x)=\int_{\partial D}\left\{\frac{\partial \Phi(x, \tilde{x})}{\partial \nu(\tilde{x})}-i \eta \Phi(x, \tilde{x})\right\} \varphi(\tilde{x}) d s(\tilde{x}), \quad x \in \mathbb{R}^{2} \backslash \overline{D},
\end{align}
where $\eta$ is the real coupling parameter and $\varphi$ the unknown density function. Then the main goal for  the direct scattering problem is to determine the unknown density function $\varphi$ according to the sound soft boundary condition by
\begin{align}\label{ac2.11}
(I+K-i \eta S) \varphi=-2 u^{\rm i} \quad \text{on} \quad \partial D.
\end{align}
By \eqref{ac2.11}, the far field pattern further can be written as
\begin{align}\label{ac2.12}
u_{\infty}(\hat{x}, d)=\frac{e^{-i \frac{\pi}{4}}}{\sqrt{8 \pi k}} \int_{\partial D}(k \nu(\tilde{x}) \cdot \hat{x}+\eta) e^{-i k \hat{x} \cdot \tilde{x}} \varphi(\tilde{x}) d s(\tilde{x}).
\end{align}
Equations \eqref{ac2.11}, \eqref{ac2.12} together with \eqref{ac2.6} give the forward operator $G$.

\section{Bayesian inversion}

Denote $X=C^{2, \alpha}[0, 2\pi)$ and $Y=\mathbb{R}^m$.  We consider the observation $y\in Y$ is a Poisson point process \cite{werner} in $[0, 2\pi)$ with parameters $\lambda_j=(G(q))_j=(\tau |u_\infty|^2)_j,\,\, j=1, 2, \cdots, m$. The goal of statistical inversion is to explore the posterior distribution $\mu^y(q)$, i.e., the distribution of $q$ conditional on the value of $y$. Bayes' formula shows that
\begin{align}
\label{biv3.1}
\frac{d\mu^{y}}{d\mu_{\rm{pr}}}=\pi(y|q),
\end{align}
where $\mu_{\rm{pr}}$ is the prior distribution for the unknown $q$ and $\pi(y|q)$ the likelihood function. 

 The Poisson noise yields the likelihood function 
\begin{align}\label{biv3.2}
\pi(y|q)=\prod_{j=1}^m \frac{\lambda_j^{y_j}e^{-\lambda_j}}{y_j!}.
\end{align}
Following \cite{stuart1, stuart2}, the likelihood function can be written in the form of 
\begin{align}
\label{biva3.3}
\pi(y|q)\propto \exp(-\Lambda(q, y)),
\end{align}
where 
\begin{align}
\label{biva3.4}
\Lambda(q, y)=\sum_{j=1}^m (G(q))_j  -y_j\log (G(q))_j.
\end{align}

The prior is given before acquiring the data. Since the $q$ denotes the length, an additional restriction should be imposed: the unknown must be positive. For this purpose, we consider the following two ways  \cite{bui, zhou}
\begin{align}
&q(t)=\exp(r(t)),\label{biv3.3}\\
&q(t)=\frac{a}{2}({\rm{erf}}(r(t))+b),\label{biv3.4}
\end{align}
where $a>0$ and $b>1$ are constants and $\rm{erf}(\cdot)$ is the error function
\begin{align*}
{\rm erf}(r(t))=\frac{2}{\sqrt{\pi}}\int_0^{r(t)}\exp(-t^2)dt.
\end{align*}
These two ways guarantee  that $q$ is strictly positive. With these parameterizations,  we now infer the new unknown parameter $r$. 
 We still use $\mu_{\rm pr}$ to denote the prior distribution for $r$. Without any confusion, the forward operator is still denoted by $G$. 
  The prior needs to be chosen carefully and it is not unique. We give different priors for \eqref{biv3.3} and \eqref{biv3.4}.

 For \eqref{biv3.3}, the prior $\mu_{\rm pr}$ is given by the Gaussian
\begin{align}\label{biv3.7}
\mu_{\rm pr}=\mu_1=N(0, C_1),
\end{align}
where $C_1$ is the covariance function and should be chosen carefully. In \cite{bui}, a convenient prior to use is 
\begin{align}\label{biv3.6}
r''(t)\sim N(0, A^{-s}), \quad s>0,
\end{align}
where $A=\frac{ d^2}{dt^2}$ with the definition domain
\begin{align}\label{biv3.7}
D(A) :=\left\{v \in H^{2}[0,2 \pi] : \quad \int_{0}^{2 \pi} v(t) d t=0\right\}.
\end{align}
In this prior assumption, it can be verified $r\in X$  almost surely with $0<\alpha<\min\{1, s-1/2\}$, i.e., $\mu_1(X)=1$ \cite{bui, stuart1}.
 In this prior assumption, the well-posedness of the posterior distribution for the Gaussian noise case has been discussed in \cite{bui}.


For \eqref{biv3.4}, the prior distribution is chosen as $\mu_{\rm pr}=N(0, C_2)$ denoted by $\mu_2$.
 To guarantee the period, we use the covariance function \cite{rasmussen, solin}
\begin{align}\label{biv3.9}
C_2(s, t)=\exp(-\frac{2\sin^2(\frac{t-s}{2})}{l^2}),
\end{align}
where $l$ defines the characteristic length-scale. The covariance example can be viewed as  mapping the one-dimensional input variable $t$  to the two-dimensional $r(t) = (\cos t, \sin t)$. The distance between $r(t)$ and $r(s)$ is given 
\begin{align*}
|r(t)-r(s)|^2=(\cos t-\cos s)^2+(\sin t-\sin s)^2=4\sin^2 \frac{t-s}{2}.
\end{align*}
The covariance kernel \eqref{biv3.9} can be obtained by taking the squared exponential kernel in $r$ space
\begin{align*}
C_2(r(t), r(s))=\exp (-\frac{|r(t)-r(s)|^2}{2l^2}).
\end{align*}
 We apply a hybrid prior   proposed in \cite{zhou}. On the basis of the Gaussian prior $\mu_2$, we further assume that the unknown $r$ satisfies the total variation prior as in \cite{zhou}
\begin{align}\label{biv3.10}
\frac{d\mu_{\rm TV}}{d\mu_2}(r)\propto \exp(-\|r\|_{\rm TV}):=\exp(-R(r)),
\end{align}
where $\|\cdot\|_{\rm TV}$ is the TV seminorm,
\begin{align}\label{biv3.11}
\|r\|_{\rm TV}=\zeta \int_0^{2\pi}|\frac{d r}{dt}|dt
\end{align}
with a positive constant $\zeta$. It follows immediately that the posterior distribution has this form
\begin{align}
\label{biv3.14}
\frac{d\mu^y}{d\mu_2}\propto\exp(-\Lambda(r, y)-R(r)):=\exp(-\Psi(r, y)).
\end{align}

Next, we discuss the well-posedness of the posterior distribution with the hybrid prior. This well-posedness is characterized in the sense of Hellinger metric, which is defined by
\begin{align*}
 d_{\rm H}(\mu^{y}, \mu^{\tilde{y}})=\left(\frac{1}{2}\int (\sqrt{\frac{d\mu^{y}}{d\mu_2}}-\sqrt{\frac{d\mu^{\tilde{y}}}{d\mu_2}})^2d\mu_2\right)^{\frac{1}{2}}.
\end{align*}
Here, we need to use the following Fernique theorem to some functionals are bounded for Gaussian measures.
\begin{thm}[Fernique]
If $\mu=N(0, C)$ is a Gaussian measure on Banach space $X$, so that $\mu(X)=1$, then there exists $\alpha> 0$ such that
\begin{align*}
\int_X \exp(\alpha\|r\|_X^2)d\mu<\infty.
\end{align*}
\end{thm}
\begin{Lemma}\label{lem3.1}
For fixed $\hat{x}, d$, there exists $M>0$ such that
\begin{align}
|G(r)| \leqslant M(1+\|r\|_X).
\end{align}
for all $r\in X$.
\end{Lemma}
\begin{Proof}
According to \eqref{ac2.6} and \eqref{ac2.12}, we have
\begin{align*}
&u_{\infty}(\hat{x})=\\
&\frac{e^{-i \frac{\pi}{4}}}{\sqrt{8 \pi k}} \int_{0}^{2 \pi}(k \nu(t) \cdot \hat{x}+\eta)e^{-i k \hat{x} \cdot\left(q \cos t, q \sin t\right)^{T}}\tilde{\varphi}(t)  \sqrt{q^2+q'^2}d t,
\end{align*}
where $\tilde{\varphi}(t)=\varphi(x)\mid_{x=x(t)}$. Hence it follows that
\begin{align*}
\left|u^{\infty}(\hat{x})\right|^2 \leqslant M_1\left\|q^2+q^{\prime 2}\right\|_{\infty} \leq M_1(\|q\|_\infty^2+\|q^{\prime}\|_\infty^2). 
\end{align*}
By \eqref{biv3.4}, we get
\begin{align*}
&\|q\|_\infty\leq \frac{a}{2}(\|{\rm erf} (r)\|_\infty+b)\leq \frac{a}{2}(1+b),\\
&\|q'\|_\infty\leq \frac{a}{\sqrt{\pi}}\|r'\exp(-r^2)\|_\infty\leq \frac{a}{\sqrt{\pi}}\|r'\|_{\infty}\leq M_2\|r\|_X.
\end{align*}
This completes the proof.
\end{Proof}
\begin{Lemma}\label{lem3.2}
For every $s>0$, there exists constants $I_2(s)\in\mathbb{R}$ and $I_1(s)>0$ such that, for all $r\in X$, and $y\in Y$ with $|y|<s$,
\begin{align}
I_2\leq\Lambda(r, y)\leq I_1.
\end{align}
\begin{Proof}
For convenient, denote 
\begin{align*}
\Lambda(r, y)&
=\langle \lambda, \mathbf{1} \rangle-\langle y, \log\lambda\rangle, 
\end{align*}
where $\langle \cdot, \cdot \rangle$ is the inner product in $\mathbb{R}^m$. It is obvious from \eqref{biv3.4} that
\begin{align*}
\frac{a}{2}(b-1)\leq q(t)\leq \frac{a}{2}(b+1).
\end{align*}
Therefore, the boundedness of the forward operator $G$ implies that
\begin{align}\label{buds}
0<\underline{\lambda}\leq \lambda\leq \overline{\lambda}
\end{align} 
for two constant vectors $\underline{\lambda}$, $\overline{\lambda}$. And it follows immediately that $|\log\lambda|\leq l_{\rm max}$ for a positive constant $l_{\rm max}$.
Furthermore, we obtain for  $|y|<s$ by the Cauchy-Schwarz inequality
\begin{align*}
&\Lambda(r, y)\geq \langle \underline{\lambda}, 1\rangle- |\log\lambda||y|\geq \langle \underline{\lambda}, 1\rangle- l_{\rm max}s:=I_2,\\
&\Lambda(r, y)\leq \langle \overline{\lambda}, 1\rangle+ |\log\lambda||y|\geq \langle \overline{\lambda}, 1\rangle+ l_{\rm max}s:=I_1.
\end{align*}
\end{Proof}
\begin{remark}
It should be noted that  the lower bound $\underline{\lambda}$ in \eqref{buds} probably has $0$ component. In this case, we can make a shift for the forward operator by
\begin{align*}
\tilde{G}(r)=G(r)+e,
\end{align*}
where $e$ is a positive constant vector. Then we discuss the new problem of the form
\begin{align*}
\tilde{G}(r)=\tilde{\lambda}:=\tau |u_\infty|^2 +e.
\end{align*}
To avoid the addition statement, we still use $G$ as the forward operator. 
\end{remark}
\end{Lemma}
\begin{thm}
Assume that $\mu_2$ is the Gaussian measure with covariance function \eqref{biv3.9}. For $y, \tilde{y}$ with $\max\{|y|, |\tilde{y}|\}<\varpi$, there exists $M=M(\varpi)>0$ such that 
\begin{align}\label{hell0}
 d_{\rm H}(\mu^{y}, \mu^{\tilde{y}})\leq M |y-\tilde{y}|.
\end{align}
\end{thm}
\begin{Proof}
First, it should be pointed out that the posterior measure $\mu^y$ is absolute continuous about the prior measure $\mu_2$. This will guarantee the existence of the Radon-Nikodym derivative $\frac{d\mu^y}{d\mu_2}$. In fact, for any measurable subset $X_1\subset X$, it holds that
\begin{align*}
\mu^y(X_1)=\int_{X_1}  \exp(-\Psi(r, y))d\mu_2\leq \int_{X_1} d\mu_2.
\end{align*}
Therefore $\mu_2(X_1)=0$ implies that $\mu^y(X_1)=0$.
Define 
\begin{align*}
A(y)=\int_X \exp(-\Psi(r, y))d\mu_2.
\end{align*}
It is obvious that 
\begin{align*}
A(y)\leq 1.
\end{align*}
By Lemma \ref{lem3.2}, we have
\begin{align*}
&A(y)\geq \int_{\|r\|_X\leq 1}\exp(-\Psi(r, y))d\mu_2\\
&=\int_{\|r\|_X\leq 1}\exp(-\Lambda(r, y)-R(r))d\mu_2\\
&=\int_{\|r\|_X\leq 1}\exp(-\sum_{j=1}^m\left\{ (G(r))_j-y_j\log(G(r))_j\right\}-\|r\|_{\rm TV})d\mu_2\\
&\geq \int_{\|r\|_X\leq 1}\exp(-I_1)\exp(-\|r\|_{\rm TV})d\mu_2\\
&\geq \exp(-I_1)\int_{\|r\|_X\leq 1}\exp(-\zeta\int_0^{2\pi} |r'(t)|dt)d\mu_2\\
&\geq \exp(-I_1)\int_{\|r\|_X\leq 1}\exp(-2\pi\zeta\|r\|_X)d\mu_2\\
&\geq \exp(-I_1)\exp(-2\pi\zeta)\int_{\|r\|_X\leq 1}d\mu_2.
\end{align*}
Since $\mu_2$ is the Gaussian measure, the unit ball has positive measure.  This shows that $A(y)$ is strictly positive. This together with this absolute continuous shows that the posterior distribution $\mu^y$ is well-defined.
The Hellinger metric holds that
\begin{align}\label{hell1}
&d^2_{\rm H}(\mu^{y}, \mu^{\tilde{y}})=\frac{1}{2}\int (\sqrt{\frac{d\mu^{y}}{d\mu_2}}-\sqrt{\frac{d\mu^{\tilde{y}}}{d\mu_2}})^2d\mu_2\\
&=\frac{1}{2}\int \left\{\left(\frac{\exp(-\Psi(r, y))}{A(y)}\right)^{\frac{1}{2}}- \left(\frac{\exp(-\Psi(r, \tilde{y}))}{A(\tilde{y})}\right)^{\frac{1}{2}}\right\}^2 d\mu_2\nonumber\\
&=\frac{1}{2}\int \left\{\left(\frac{\exp(-\Psi(r, y))}{A(y)}\right)^{\frac{1}{2}}- \left(\frac{\exp(-\Psi(r, \tilde{y}))}{A(y)}\right)^{\frac{1}{2}}\right\}^2 d\mu_2\nonumber\\
&+\frac{1}{2}\int \left\{\left(\frac{\exp(-\Psi(r, \tilde{y}))}{A(y)}\right)^{\frac{1}{2}}- \left(\frac{\exp(-\Psi(r, \tilde{y}))}{A(\tilde{y})}\right)^{\frac{1}{2}}\right\}^2 d\mu_2\nonumber\\
&\leq \frac{1}{2A(y)} \int \left\{\exp(-\frac{\Psi(r, y)}{2})-\exp(-\frac{\Psi(r, \tilde{y})}{2})\right\}^2d\mu_2\nonumber\\
&+ \frac{1}{2} |A(y)^{-\frac{1}{2}}-A(\tilde{y})^{-\frac{1}{2}}|^2\int \exp(-\Psi(r, \tilde{y}))d\mu_2:=\Sigma_1+\Sigma_2.\nonumber
\end{align}
For $\Sigma_1$, since $A(y)$ is bounded, it suffices to estimate the integral part. 
According to Lemma \ref{lem3.1} and Fernique therorem, we have
\begin{align}
&\int \left\{\exp(-\frac{\Psi(r, y)}{2})-\exp(-\frac{\Psi(r, \tilde{y})}{2})\right\}^2d\mu_2\nonumber\\
&\leq \frac{1}{4}\int |\Psi(r, y)-\Psi(r, \tilde{y})|^2d\mu_2\nonumber\\
&=\frac{1}{4}\int |\sum_{j=1}^m (\tilde{y}_i-y_i)(\log G(r))_i|^2 d\mu_2\nonumber\\
&\leq \frac{1}{4}|y-\tilde{y}|^2\int |\log G(r)|^2d\mu_2 \nonumber
\\
&\leq \frac{M^2}{4}|y-\tilde{y}|^2\int\log^2 (1+\|r\|_X) d\mu_2\leq \tilde{M} |y-\tilde{y}|^2.\label{hell2}
\end{align}
Since $A(y)$ and $A(\tilde{y})$ are bounded from below, it follows that
\begin{align}
&|A(y)^{-\frac{1}{2}}-A(\tilde{y})^{-\frac{1}{2}}|^2\nonumber\\
&\leq M \max(A(y)^{-3}, A(\tilde{y})^{-3})|A(y)-A(\tilde{y})|^2
\leq M|y-\tilde{y}|^2.\label{hell3}
\end{align}
The estimation \eqref{hell0} can be obtained from \eqref{hell1}, \eqref{hell2} and \eqref{hell3}.
\end{Proof}
\begin{remark}
The requirement $\mu_2(X)=1$ in Fernique theorm is guaranteed naturally since it has the  exponential covariance function \eqref{biv3.9}.
This case for $q$ with the form \eqref{biv3.3} is similar to \eqref{biv3.4}. Actually, we just need to take $\zeta=0$ and verify the condition $\mu_1(X)=1$. 
\end{remark}

\section{Numerical test}

In this section, we demonstrate some numerical examples to show the effectiveness of the Bayesian method for the inverse obstacle scattering problem with the Poisson data. We are ready to use the MCMC method to generate samples to explore the posterior distribution $\mu^y$ in \eqref{biv3.1}. 

It is clear that the operator $A$ in \eqref{biv3.7} has the eigen-system $\{n^2; \frac{\cos(nt)}{\sqrt{\pi}}, \frac{\sin(nt)}{\sqrt{\pi}}\}$. With this, the prior for \eqref{biv3.3} can be generated by the Karhunen-Lo\`eve expansion. 
 The truncated expansion is used to generate the prior samples. The corresponding derivatives can be computed accurately. 
For \eqref{biv3.4}, the samples for $r$ with the covariance \eqref{biv3.9} are period. We use the fast Fourier transform (FFT) to implement the numerical computations for the derivatives of the samples. The posterior samples for the two ways are both generated by the Metropolis-Hastings algorithm. 

%

We just list the hybrid prior case in the following. For using \eqref{biv3.3} and the corresponding prior, the process is similar. 
\begin{itemize}
\item Initialize $r^0$ and further generate $q^0$ in \eqref{biv3.4}. 
\item For iteration $i=1, 2, \cdots$ do
\begin{itemize}
\item Propose $\hat{r}=\sqrt{1-\beta^2} r^{i-1}+\beta \xi, \,\, \xi\sim N(0, C_2)$ and further generate $\hat{q}$.
\item Acceptance probability: 
\begin{align*}
&\alpha(\hat{r}; r^{i-1})=\min\left\{1, \frac{\exp(-\Psi(\hat{r}, y))}{\exp(-\Psi(r^{i-1}, y))}\right\},\\
&\upsilon\sim {\rm Uniform}(0, 1).
\end{align*}
\item Accept or reject the proposal in the following way
\begin{align*}
r^{i}=\left\{
\begin{aligned}
&\hat{r}, \,\,\,\,\,\,\,\,\,\,{\rm if}\,\,  \upsilon>\alpha,\\
&r^{i-1}, \,\, {\rm else}.
\end{aligned}
\right.
\end{align*}
\end{itemize}
endfor
\end{itemize}
Consider the following obstacles
\begin{align*}
&1: {\rm Peanut:}\,\,  q(t)=3(\cos^2 t+0.25\sin^2 t), x_1(t)=q(t)\cos t, x_2(t)=q(t)\sin t;\\
&2: {\rm Kite:}\,\, \left(x_{1}, x_{2}\right)=(\cos t+0.65 \cos 2 t-0.65, 1.5 \sin t);\\
&3: {\rm Drop:} \,\,x_1=-1+2\sin(t/2), x_2=-\sin(t);\\
&4: {\rm Cloverleaf:} \,\,q(t)=1+0.3\cos(4t), x_1(t)=q(t)\cos t, x_2(t)=q(t)\sin(t),
\end{align*}
where  $t \in[0,2 \pi)$.

\begin{figure*}[!t]
\centering
\subfigure[Penut data] {\includegraphics[height=1.5in,width=2in,angle=0]{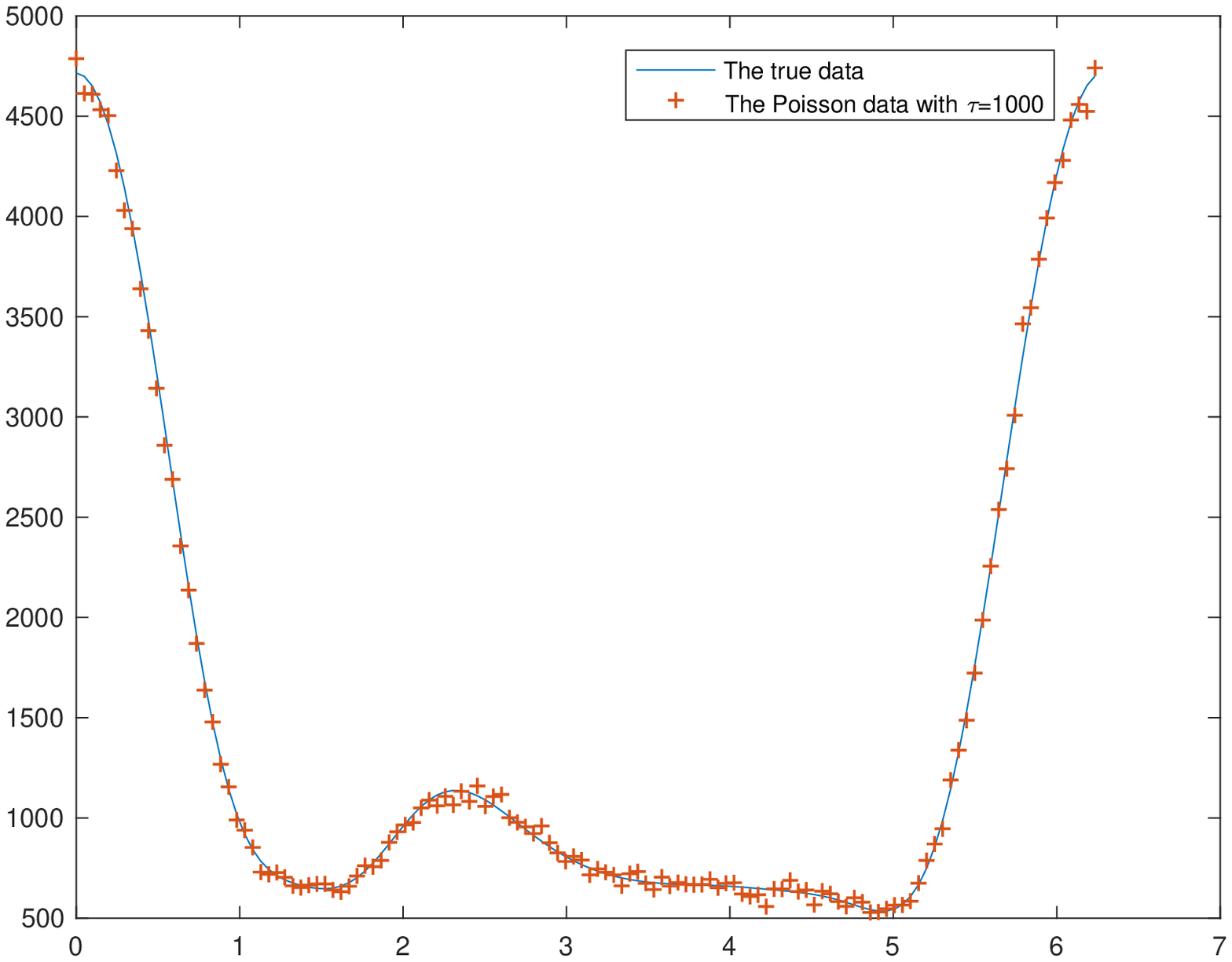}}
\subfigure[Kite data] {\includegraphics[height=1.5in,width=2in,angle=0]{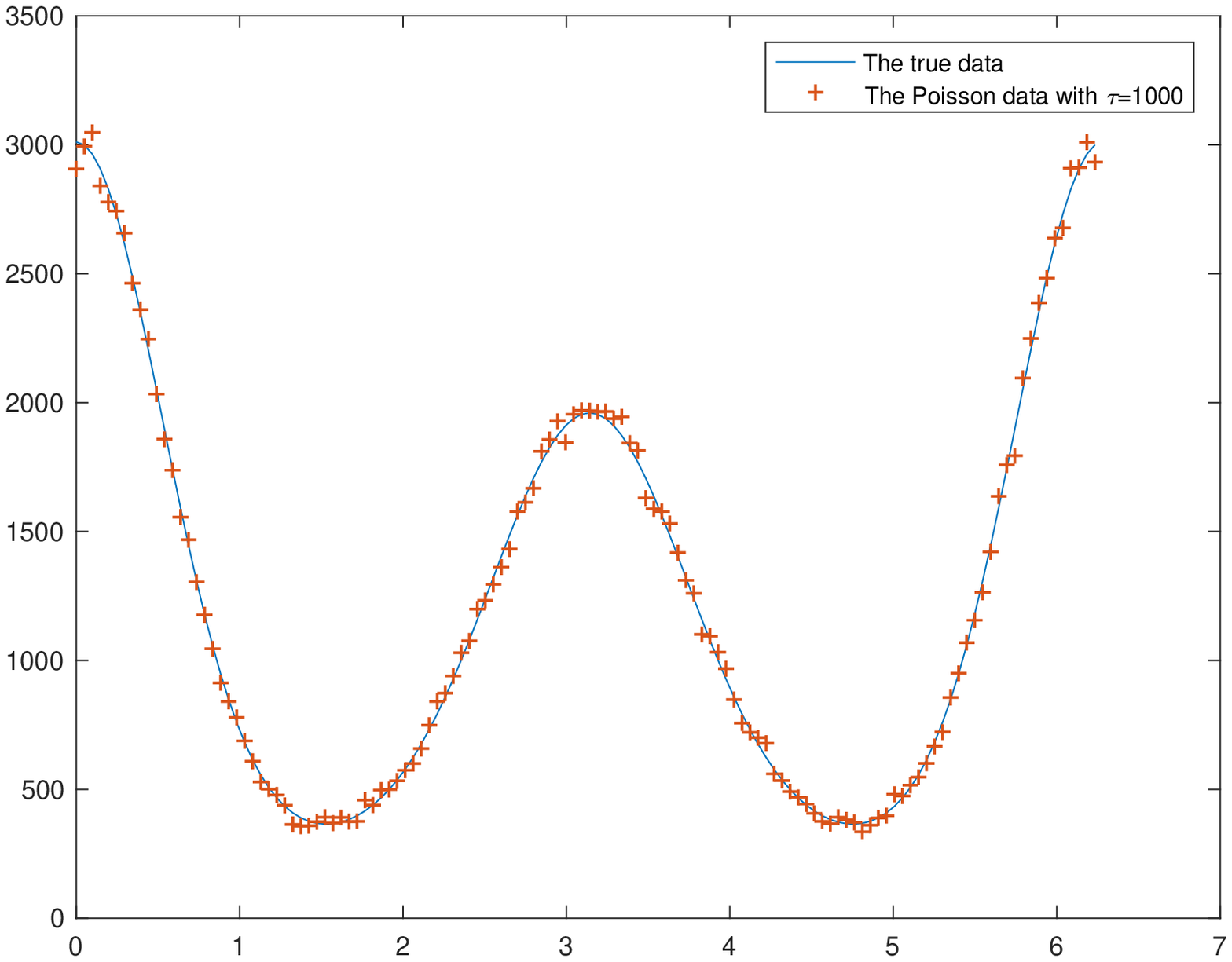}}
\subfigure[Cloverleaf data] {\includegraphics[height=1.5in,width=2in,angle=0]{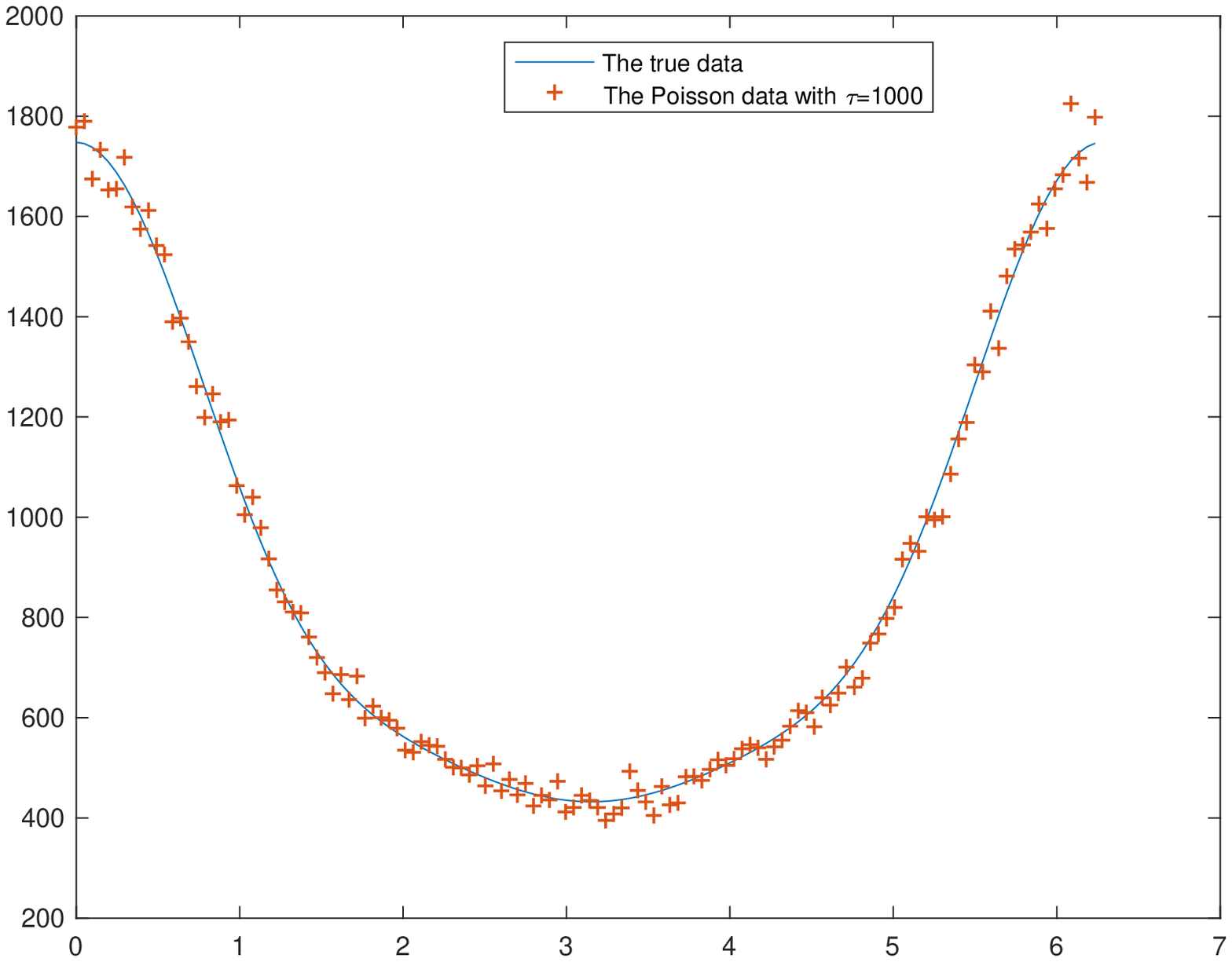}}
\subfigure[Drop data] {\includegraphics[height=1.5in,width=2in,angle=0]{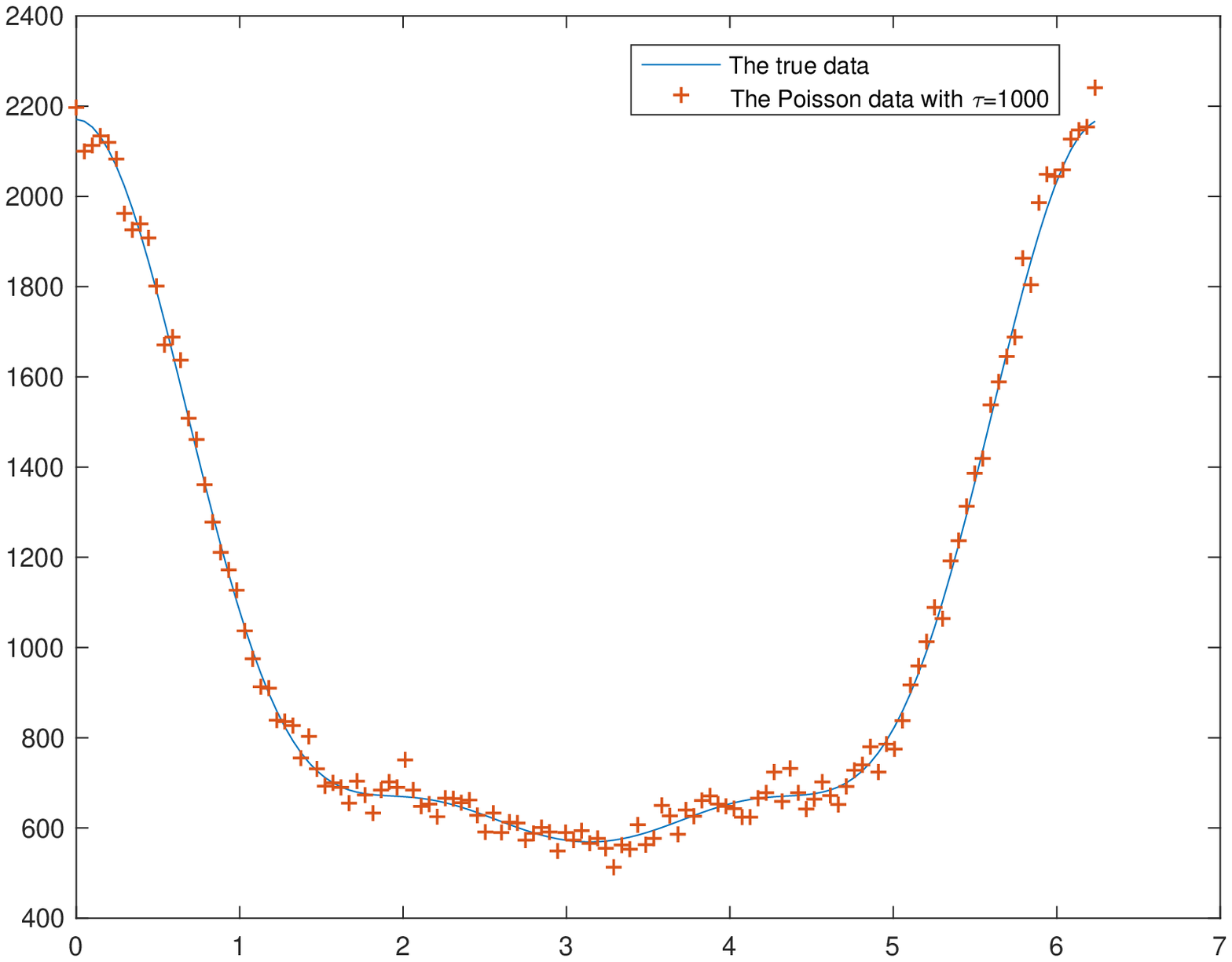}}
\caption{The Poisson data with $\tau=1000$.}
\label{data1}
\end{figure*}

\begin{figure*}[!t]
\centering
\subfigure[Gaussian prior] {\includegraphics[height=1.5in,width=2in,angle=0]{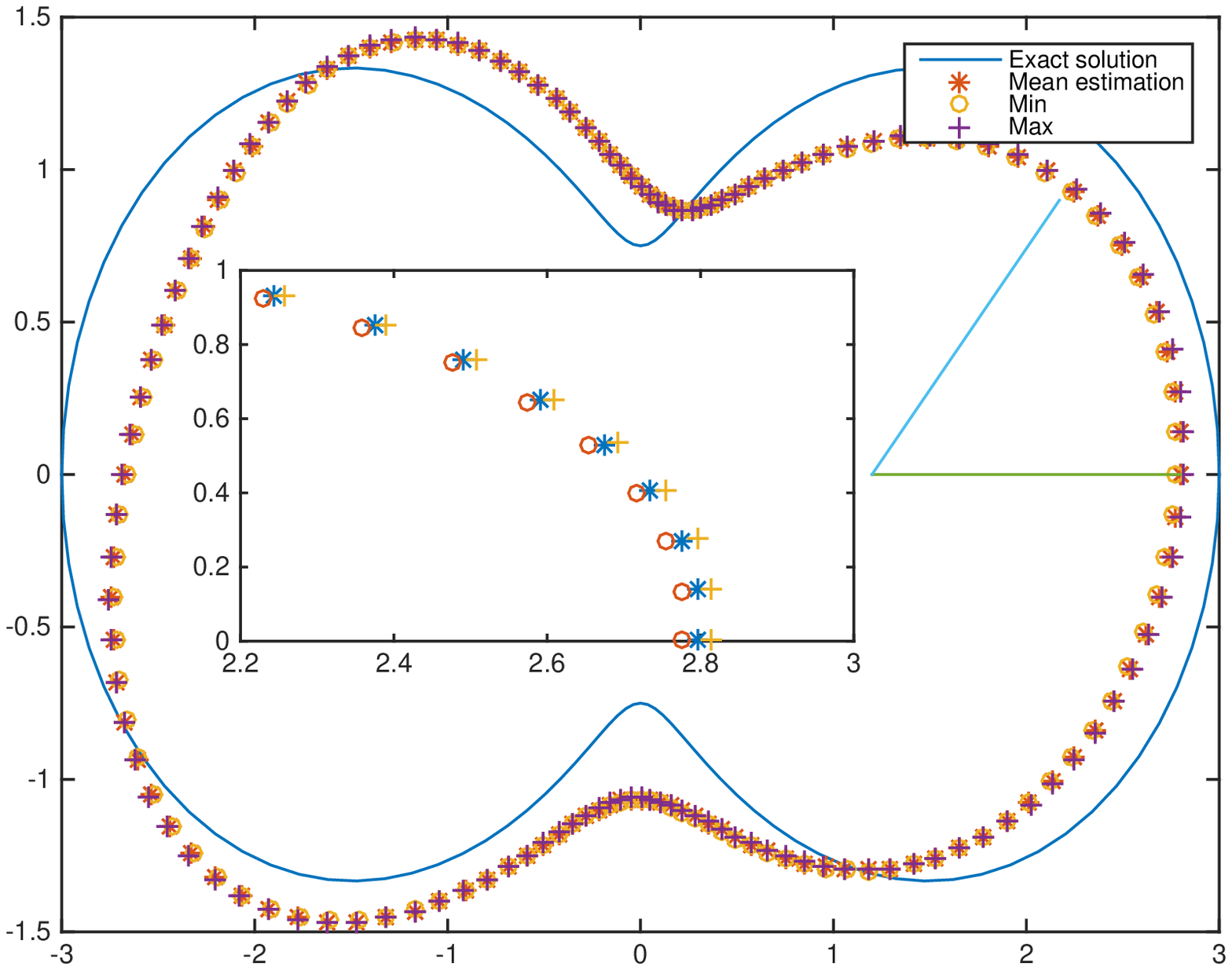}}
\subfigure[Hybrid total variation prior] {\includegraphics[height=1.5in,width=2in,angle=0]{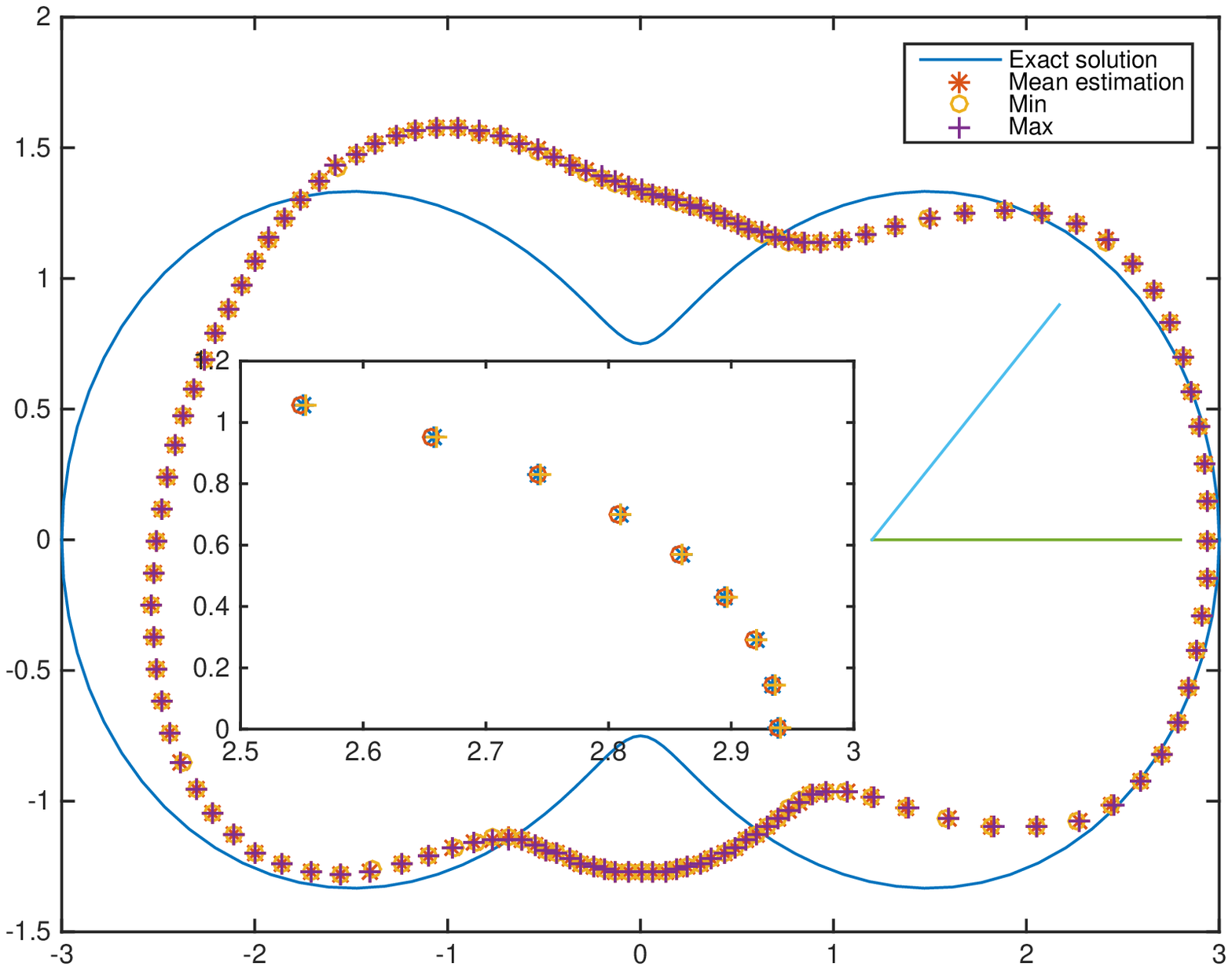}}
\caption{The numerical reconstructions of the peanut shape}
\label{fig1}
\end{figure*}

\begin{figure*}[!t]
\centering
\subfigure[Drop] {\includegraphics[height=1.5in,width=2in,angle=0]{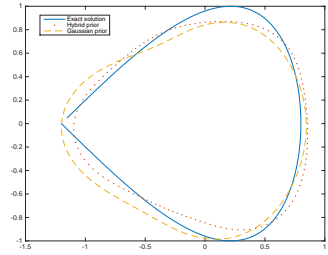}}\\
\subfigure[Kite] {\includegraphics[height=1.5in,width=2in,angle=0]{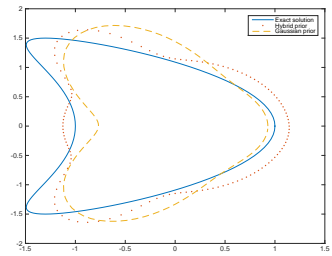}}\\
\subfigure[Cloverleaf] {\includegraphics[height=1.5in,width=2in,angle=0]{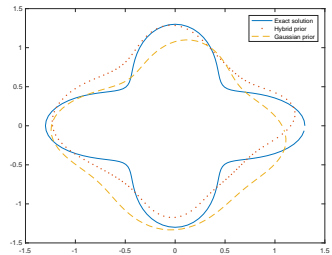}}
\caption{The numerical reconstructions of the other three obstacles}
\label{fig2}
\end{figure*}

In Fig. \ref{data1}, we plot the Poisson data for $\tau=1000$.  Using these data, we reconstruct the obstacles  displayed  in Fig. \ref{fig1} and \ref{fig2}.
The posterior sample means are used as the approximations. In Fig. \ref{fig1}, we plot the $95\%$ confidence interval for the peanut shape.

\end{document}